\documentclass{amsart}

\usepackage{amssymb}
\usepackage{amsthm}
\usepackage{amsmath}
\usepackage{amsbsy}

\usepackage{tikz}
\newtheorem*{theorem}{Theorem}
\newtheorem*{proposition}{Proposition}

\newtheorem*{definition}{Definition}

\newcommand{\Abs}[1]{\left\lvert#1\right\rvert}

\begin{document}

\title[]{Synchronization of Kuramoto oscillators\\ in Dense Networks}

\thanks{J.L. was supported in part by the National Science Foundation under award DMS-1454939 and CCF-1934964. S.S. was partially supported by the NSF (DMS-1763179) and the Alfred P. Sloan Foundation. }

\author[]{Jianfeng Lu}
\address[Jianfeng Lu]{Department of Mathematics, Department of Physics, and Department of Chemistry,
Duke University, Box 90320, Durham NC 27708, USA}
\email{jianfeng@math.duke.edu}

\author[]{Stefan Steinerberger}
\address[Stefan Steinerberger]{Department of Mathematics, Yale University, New Haven, CT 06510, USA}
\email{stefan.steinerberger@yale.edu}

\begin{abstract} We study synchronization properties of systems of  Kuramoto oscillators. The problem can also be understood as a question about
the properties of an energy landscape created by a graph. More formally, let $G=(V,E)$ be a connected graph and $(a_{ij})_{i,j=1}^{n}$ denotes its adjacency matrix. Let the function $f:\mathbb{T}^n \rightarrow \mathbb{R}$ be given by
$$ f(\theta_1, \dots, \theta_n) = \sum_{i,j=1}^{n}{ a_{ij} \cos{(\theta_i - \theta_j)}}.$$
This function has a global maximum when $\theta_i = \theta$ for all
$1\leq i \leq n$. It is known that if every vertex is connected to at
least $\mu(n-1)$ other vertices for $\mu$ sufficiently large, then
every local maximum is global. Taylor proved this for
$\mu \geq 0.9395$ and Ling, Xu \& Bandeira improved this to
$\mu \geq 0.7929$. We give a slight
improvement to $\mu \geq 0.7889$.  Townsend, Stillman \& Strogatz suggested that the critical value might be $\mu_c = 0.75$.
\end{abstract}
\maketitle

\section{Introduction}
We study a simple problem that can be understood from a variety of
perspectives. Perhaps its simplest formulation is as follows: let
$G=(V,E)$ be a connected graph and $(a_{ij})_{i,j=1}^{n}$ denotes its
adjacency matrix. We assume the graph is simple, and thus
  $a_{ii} = 0$ for $i = 1, \ldots, n$. We are then interested in the
behavior of the energy functional
$f:\mathbb{T}^n \cong [0,2\pi]^n \rightarrow \mathbb{R}$ given by
\begin{equation} \label{main}
 f(\theta_1, \dots, \theta_n) = \sum_{i,j=1}^{n}{ a_{ij} \cos{(\theta_i - \theta_j)}}.
 \end{equation}
Ling, Xu \& Bandeira \cite{ling} ask the following very interesting
\begin{quote}
  \textbf{Question.} What is the relationship between the existence of
  local maxima and the topology of the network?
\end{quote}

$f$ assumes its global maximum when $\theta_i \equiv \theta$ is constant and this is the unique global maximum up to rotation. Factoring out the rotation symmetry, there are at least $2^n$ critical points of the form $\theta_i \in \left\{0, \pi\right\}$. The main question is under which condition we can exclude the existence of local maxima that are not global maxima.\\
This is related to the Kuramoto model as follows: suppose we consider the system of ordinary differential equations given by
$$ \frac{d \theta_i}{dt} = -\sum_{j=1}^{n}{a_{ij} \sin{(\theta_i - \theta_j)}}.$$
We can interpret this system of ODEs as a gradient flow with respect to the energy
$$ E(\theta_1, \dots, \theta_n) =  - \sum_{i,j=1}^{n}{ a_{ij} \cos{(\theta_i - \theta_j)}}.$$
In this case, local maxima that are not global correspond to stable local minima of the gradient flow. In light of this model, particles on the circle that are connected by springs, it is natural to assume that no spurious local minimizer of this energy exist if there are enough springs. This motivated an existing line of research: Taylor \cite{taylor} proved that if each vertex is connected to at least $\mu(n-1)$ vertices for $\mu \geq 0.9395$, then \eqref{main} does not have local maxima that are not global. Ling, Xu \& Bandeira \cite{ling} improved this to $\mu \geq 0.7929$. They also showed the existence of a configuration coming from the family of Wiley-Strogatz-Girvan networks \cite{wiley} where each vertex is connected to $0.68n$ other vertices that indeed has local maxima that are not global. Townsend, Stillman \& Strogatz \cite{townsend} suggest that the critical value might be $\mu_c = 0.75$ and identify networks with $\mu = 0.75$ having interesting spectral properties. \\

The problem itself arises in a variety of settings. We refer to the surveys \cite{review1, review2, review3} for an overview regarding synchronization problems, to \cite{chaos} for insights into complexities of the Kuramoto model and to \cite{lopes, lopes2} for random Kuramoto models. There is also recent interest in the landscape of non-convex loss functionals for which this problem is a natural test case, we refer to \cite{loss1, loss2, loss3, loss4}.

\begin{theorem} If $G=(V,E)$ is a connected graph such that the degree of every vertex is at least $0.7889(n-1)$,  then
$$ f(\theta_1, \dots, \theta_n) = \sum_{i,j=1}^{n}{ a_{ij} \cos{(\theta_i - \theta_j)}}$$
does not have local maxima that are not global.
\end{theorem}

The main idea behind the argument is a refinement of the approach of
Ling, Xu \& Bandeira \cite{ling} in a certain parameter range using a
new decomposition of the points. We consider the problem a natural
benchmark for testing our understanding of the geometry of energy
landscapes. We conclude by reiterating the original question from
\cite{ling}: which kind of assumption on the network (this paper, for
example, is only dealing with edge-density assumptions) implies
synchronization?

\section{Proof}
\subsection{Ingredients.} The purpose of this section is to sketch several of the tools that go into the argument which is a variation on the argument given by Ling, Xu \& Bandeira \cite{ling}. We first recall their argument. They start by introducing a useful Proposition (precursors of which can be found in Taylor \cite{taylor}).

\begin{proposition}[Ling, Xu, Bandeira \cite{ling}] Let $(\theta_1, \dots, \theta_n) \in \mathbb{T}^n$ be a strict local maximizer of \eqref{main}. If there exists an angle $\theta_r$ such that 
$$ \forall~1 \leq i \leq n: \qquad |\sin{(\theta_i - \theta_r)}| < \frac{1}{\sqrt{2}},$$
then all the $\theta_i$ have the same value.
\end{proposition}
The idea behind this argument is as follows: if
$(\theta_1, \dots, \theta_n) \in \mathbb{T}^n$ is a local maximizer,
then the quadratic form (corresponding to the negative Hessian) is positive semi-definite.
In other words, a necessary condition for being a local
maximum (derived in \cite{ling}) is that for all vectors
$w \in \mathbb{R}^n$
$$ \sum_{i,j=1}^{n}{ a_{ij} \cos{(\theta_i - \theta_j)} (w_i - w_j)^2} \geq 0.$$
We can derive a contradiction by defining a vector
$w \in \left\{-1, 1\right\}^n$ depending on which of the two `cones'
the variable $\theta_i$ is in. 
Then the summation only ranges over pairs that are in opposite sides of the cone. The cosine is negative for those values and since the graph is connected, there
is at least one connection leading to a contradiction.\\

A second important ingredient is the Kuramoto parameter \cite{kuramoto}
$$ r = \| r\| e^{i \theta_r}:= \sum_{j=1}^{n}{ e^{i \theta_j}}.$$
The second part in the argument \cite{ling} is based on showing that
\begin{equation} \label{length}
\begin{aligned}
\|r\|^2 &\geq \frac{n^2}{2} - \sum_{i \neq j}{ (1-a_{ij}) \Abs{ \cos{(\theta_i - \theta_j)} - \cos^2{(\theta_i - \theta_j)} }}  \\
&\geq \left(2 \mu - \frac{3}{2} \right)n^2 + 2(1-\mu)n,
\end{aligned}
\end{equation}
where the second inequality follows from
$$ \Abs{\cos{(\theta_i - \theta_j)} - \cos^2{(\theta_i - \theta_j)}} \leq 2$$
and 
$$ \sum_{i \neq j}{(1-a_{ij})} = \sum_{i=1}^{n} \sum_{j \neq i}{(1-a_{ij})} \leq \sum_{i=1}^{n}{ (1-\mu)(n-1)} = (1-\mu) (n-1)n.$$

The argument in \cite{ling} proceeds by writing
$$ \| r\| e^{-i (\theta_i - \theta_r)} = r e^{- i\theta_i} = \sum_{j=1}^{n}{ e^{-i(\theta_i - \theta_j)}}$$
and taking imaginary parts to obtain
$$ \|r\| \sin{(\theta_i - \theta_r)} = \sum_{j=1}^{n}{ \sin{(\theta_i - \theta_j)}}.$$
However, the first order condition in a maximum implies
$$ \sum_{j=1}^{n}{a_{ij} \sin{(\theta_i - \theta_j)}} = 0$$
and thus we obtain
$$ \|r\| \sin{(\theta_i - \theta_r)} = \sum_{j=1}^{n}{(1-a_{ij}) \sin{(\theta_i - \theta_j)}}.$$
As a consequence, we have
\begin{equation} \label{bound}
 | \sin{(\theta_i - \theta_r)} | \leq \frac{1}{\|r\|} \Biggl\lvert \sum_{j=1}^{n}{(1-a_{ij}) \sin{(\theta_i - \theta_j)}} \Biggr\rvert \leq \frac{(1-\mu)n}{\|r\|}.
\end{equation}
The Proposition together with \eqref{length} and \eqref{bound} then imply the result.

\subsection{The Proof} 
Our proof is motivated by the following simple observation: inequality \eqref{length} is only sharp when 
$$ \Abs{\cos{(\theta_i - \theta_j)} - \cos^2{(\theta_i - \theta_j)}} = 2$$
for all pairs of angles $\theta_i, \theta_j$ for which $a_{ij} = 0$.  This only happens when $\theta_i - \theta_j = \pi$.
So the only extreme case of equality is when the points that are not connected are concentrated on two different sides of the torus. 
Then, however, we can dramatically improve inequality \eqref{bound} since $\sin{(\pi)} = 0$. 
The proof will decouple into two steps: either inequality \eqref{length}
is very far from sharp, in that case the origin Ling-Xu-Bandeira argument will result in the desired improvement, or
the inequality is close to being sharp in which case we can try to extract more information from it.\\

Our proof makes use of several new parameters. 
As will come to no
surprise to the reader, we obtain them by working with unspecified
coefficients in the beginning and then solving the arising
optimization problem to obtain the optimal selection of
parameters. For readers who prefer explicit values to have an idea of
scales, we will later set
$$ \varepsilon = 0.5 \qquad \mbox{and} \qquad \delta = 0.88.$$

\begin{proof} The proof decouples into several steps.\\

\textbf{Step 1: Finding Good Vertices.}
 The first part of our proof emulates the argument from \cite{ling} with one slight modification. We assume that we are working with a slightly improved bound in \eqref{length}. Specifically, we first assume that we have the identity 
$$
\sum_{i \neq j} (1 - a_{ij}) \Abs{\cos(\theta_i - \theta_j) -
\cos(\theta_i-\theta_j)^2} = (2-\alpha) (1-\mu) (n-1)n
$$
for some value of $\alpha > 0.0537$. In that case, running the original Ling-Xu-Bandeira argument again shows that we have 
$$ 
| \sin{(\theta_i - \theta_r)} | < \frac{1}{\sqrt{2}}$$ as long as
$\mu \geq 0.788897$. It thus remains to obtain a similar bound in the
case where this assumption does not hold. We may thus additionally
assume
\begin{equation} \label{standingassumption}
\sum_{i \neq j} (1 - a_{ij}) \Abs{\cos(\theta_i - \theta_j) -
\cos(\theta_i-\theta_j)^2} = (2-\alpha) (1-\mu) (n-1)n,
\end{equation}
where $\alpha$ is defined by the equation and satisfies $0 \leq \alpha \leq 0.0537$. The first assumption that we derive from this assumption is that there are many `good' vertices, a term that we will use repeatedly and that we now formally define.
\begin{definition} A vertex $i$ is `$\varepsilon-$good' if
\begin{equation} \label{firstcond}
\sum_{j=1,\, j \neq i}^{n} (1 - a_{ij})\Abs{\cos(\theta_i - \theta_j) -
\cos(\theta_i-\theta_j)^2} \geq (2-\varepsilon) (1-\mu)(n-1).
\end{equation}
\end{definition}
In practice, we will assume that $\varepsilon$ is fixed and to be optimized over later and we will refer to them as `good vertices' or `good points'. The nomenclature is motivated by the fact that the double sum has to be large (see equation \eqref{standingassumption}) and `good points' contribute to achieving this goal. It will also be reflected in the fact that we will derive some improved bounds for good points early on in the paper.

The first step in our argument is to conclude that, for every given $\varepsilon \in [\alpha, 2]$, there are at least $(1-\frac{\alpha}{\varepsilon})n$ vertices that are $\varepsilon-$good.
Suppose this is false, then the double sum in \eqref{standingassumption}
could be bounded from above by
\begin{align*}
\text{LHS of }\eqref{standingassumption}&<  \frac{\alpha}{\varepsilon} n (2-\varepsilon) (1-\mu)(n-1) + \Bigl(1-\frac{\alpha}{\varepsilon}\Bigr)n 2 (1-\mu)(n-1) \\
 &= (1-\mu)(n-1)n \Bigl( \frac{\alpha}{\varepsilon}(2-\varepsilon) + 2-2\frac{\alpha}{\varepsilon}\Bigr)\\
 &= (2-\alpha)(1-\mu)(n-1)n
\end{align*}
which is a contradiction to \eqref{standingassumption}.  \\

Let us now take a $\varepsilon-$good vertex $i$. We argue that, for
every given $ \delta \in (\varepsilon, 2)$ there are many
non-neighbors, indices $j$ such that $a_{ij} = 0$, for which
$$ \Abs{\cos(\theta_i - \theta_j) -
\cos(\theta_i-\theta_j)^2}  \geq 2 - \delta.$$
Let us assume their number is $(1-\mu - c)(n-1)$. Then we can bound, using the fact that the total number of non-neighbors is at most $(1-\mu)(n-1)$, that
$$ \sum_{j=1}^{n} (1 - a_{ij})\Abs{\cos(\theta_i - \theta_j) -
\cos(\theta_i-\theta_j)^2} \leq (1-\mu-c)(n-1) 2 + c(n-1)(2-\delta).$$
We require, using \eqref{firstcond}, that 
$$  (2-\varepsilon) (1-\mu)\leq (1-\mu-c) 2 + (2-\delta)c$$
and thus
$$  c \leq (1-\mu)\frac{\varepsilon}{\delta}.$$
This shows that for any $\varepsilon-$good vertex, the number of non-neighbors for which the cosine quantity exceeds $2-\delta$ is at least
$$(1-\mu - c)(n-1)  \geq (1-\mu)(n-1)\left(1- \frac{\varepsilon}{\delta} \right).$$

We summarize our arguments up to this step.
\begin{enumerate}
\item Let us consider the value $\alpha$ as defined in \eqref{standingassumption}. If $\alpha > 0.0537$, then we get the desired result directly from the argument of Ling, Xu \& Bandeira. It thus remains to study the cases where $0 \leq \alpha \leq 0.0537$.
\item In this case, for each $\varepsilon \geq \alpha$, there are at least $(1-\frac{\alpha}{\varepsilon})n$ vertices $i$ (`the $\varepsilon$-good vertices') for which we have the inequality
$$
\sum_{j=1,j\neq i}^{n} (1 - a_{ij})\Abs{\cos(\theta_i - \theta_j) -
\cos(\theta_i-\theta_j)^2} \geq (2-\varepsilon) (1-\mu)(n-1).
$$
\item For each $\varepsilon < \delta < 2$, each of these $(1-\frac{\alpha}{\varepsilon})n$ good points has at least $(1-\mu) (n-1)(1 - \frac{\varepsilon}{\delta})$ non-neighbors, $a_{ij} = 0$, for which
$$ \Abs{\cos(\theta_i - \theta_j) - \cos(\theta_i-\theta_j)^2}  \geq 2- \delta.$$
\end{enumerate}

\textbf{Step 2: Improved Bounds for Good Vertices.} It is an
elementary trigonometric fact that if $0 \leq \delta \leq 1$ (in
  fact the inequality below holds for $0 \leq \delta < 7/4$), then
$$ \Abs{\cos(x) - \cos(x)^2}  \geq  2 - \delta,\quad \text{implies} \quad \left|\sin{(x)}\right| \leq \frac{1}{\sqrt{2}} \sqrt{ \sqrt{9 - 4 \delta} - 3 + 2\delta} =: s_{\delta},$$
where we introduced the shorthand $s_{\delta}$ for simplicity of exposition. 
Combining these facts, we can show that for all the good points, of which there are at least $(1-\frac{\alpha}{\varepsilon})n$, we have
\begin{align*}
\|r\| \cdot \left|\sin{(\theta_i - \theta_r)}\right| &= \left| \sum_{j=1}^{n}{ (1-a_{ij}) \sin{(\theta_i - \theta_j)}} \right|\\
&\leq (1-\mu)(n-1) \Bigl(1-\frac{\varepsilon}{\delta}\Bigr) s_{\delta}  \\
&\qquad + (1-\mu)(n-1)\frac{\varepsilon}{\delta}\\
&=  (1-\mu)(n-1)\Bigl( s_{\delta} + \frac{\varepsilon}{\delta} (1-s_{\delta})\Bigr).
\end{align*}
This, in turn, implies that any good point $i$ satisfies
$$  \left|\sin{(\theta_i - \theta_r)}\right|^2 \leq \frac{ \left( s_{\delta} + \frac{\varepsilon}{\delta} (1-s_{\delta})\right)^2 (1-\mu)^2 n^2}{\|r\|^2}$$
Note that we also have \eqref{length} and \eqref{standingassumption} implying
\begin{align*}
 \|r\|^2 \geq  \left(\frac{1}{2} - (2-\alpha)(1-\mu) \right)n^2
 \end{align*}
Therefore, if $i$ is a good point, then
$$  \left|\sin{(\theta_i - \theta_r)}\right|^2 \leq \frac{ \left( s_{\delta} + \frac{\varepsilon}{\delta} (1-s_{\delta})\right)^2 (1-\mu)^2}{  \left( \frac{1}{2} - (2-\alpha)(1-\mu)\right)}.$$
Recall also that there are at most $\frac{\alpha}{\varepsilon} n$ `bad' points (points which are not good, we will also call then `outliers'). 
We define 
$$\varphi = \varphi(\mu, \varepsilon, \delta, \alpha)$$ 
as the positive angle satisfying
$$  \left|\sin{(\varphi)}\right|^2 = \frac{ \left( s_{\delta} + \frac{\varepsilon}{\delta} (1-s_{\delta})\right)^2 (1-\mu)^2}{  \left( \frac{1}{2} - (2-\alpha)(1-\mu)\right)}.$$
This is the bound we get on the angle that good points have with $\theta_r$. There is a balancing act: by choosing the parameters in a more restricted fashion, we can get a better upper bound for this angle but there will be less good points. Choosing the parameters in the other direction will result in a bigger number of good points but less control on their geometric distribution.
We will require, further along in the argument, that $\varepsilon$ and $\delta$ are chosen in such a way that
\begin{equation} \label{sincond}
  \left|\sin{(\varphi)}\right|^2 = \frac{ \left( s_{\delta} + \frac{\varepsilon}{\delta} (1-s_{\delta})\right)^2 (1-\mu)^2}{  \left( \frac{1}{2} - (2-\alpha)(1-\mu)\right)} < \frac{1}{2}.
  \end{equation}

\textbf{Step 3: The Distribution of Good Points.}
We now introduce a bit of orientation and assume, without loss of generality after possibly rotating all the points, that
$$ r = \| r\| e^{i \theta_r} = \sum_{j=1}^{n}{ e^{i \theta_j}} \qquad \mbox{is a positive real number.}$$
We know that the good points, of which there are at least $(1-\frac{\alpha}{\varepsilon})n$ many, satisfy \eqref{sincond} and are thus contained in one of two cones. 
We assume that we have $\gamma_1 n$ good points in the left cone and $\gamma_2n$ good points in the right cone (see Figure~\ref{fig:cones} for a sketch of this). Without loss of generality, we can reflect the picture if necessary and assume that $\gamma_2 \geq \gamma_1$. Recalling our lower bound on the number of good points, we have
\begin{equation}\label{eq:gamma12}
  \gamma_1 + \gamma_2 \geq 1 - \frac{\alpha}{\varepsilon}.
  \end{equation}

\begin{center}
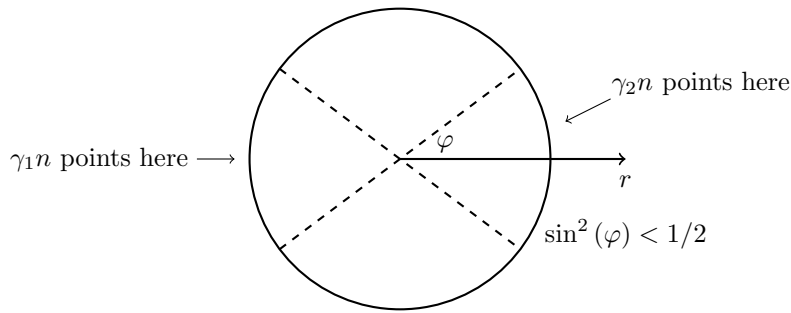
\begin{figure}[h!]
\begin{tikzpicture}
\draw [thick] (0,0) circle (2cm);
\draw [->, thick] (0,0) -- (3,0);
\node at (3, -0.3) {$r$};
\draw [thick, dashed] (-1.6, -1.2) -- (1.6, 1.2);
\draw [thick, dashed] (-1.6, 1.2) -- (1.6, -1.2);
\node at (-4,0) {$\gamma_1 n$ points here};
\draw[->] (-2.7,0) -- (-2.2,0);
\node at (4,1) {$\gamma_2 n$ points here};
\draw [->] (2.8, 0.8) -- (2.2, 0.5);
\node at (0.6, 0.2) {$\varphi$};
\node at (3,-1) {$\sin^2{(\varphi)} < 1/2$};
\end{tikzpicture}
\caption{Introducing orientation: $r$ being a positive real forces all the good points to be in two cones. The outliers can be anywhere (and could also be in the cone).\label{fig:cones}}
\end{figure}
\end{center}

We note that the outliers, of which there are (at most) $\frac{\alpha}{\varepsilon} n$, might also be in the left or the right cone, we do not make any statement about their actual location and will always assume
that they are working against us.
The inequality
$$ \|r\|^2 \geq  \left(\frac{1}{2} - (2-\alpha)(1-\mu) \right)n^2$$
forces some restrictions on $\gamma_1$ and $\gamma_2$: in particular, if all the good points and all the outliers were distributed somewhat evenly, then $\|r\|$ would actually be quite small. However, we do have a lower bound on $\|r\|$ and this forces some restrictions which we will now explore.
 Assuming the worst case (where all the outliers are actually working in our favor and contribute to making $r$ as big as possible), we have
$$ \|r\| \leq \left(\gamma_2 - \cos{(\varphi)} \gamma_1 + \frac{\alpha}{\varepsilon}\right) n$$
and therefore
$$ \gamma_2 - \cos{(\varphi)} \gamma_1 + \frac{\alpha}{\varepsilon} \geq  \left( \frac{1}{2} - (2-\alpha)(1-\mu)\right)^{1/2}$$
which implies, using $-\gamma_1 \leq \gamma_2 -1 + \frac{\alpha}{\varepsilon}$,
$$ \bigl(1 + \cos(\varphi)\bigr) \Bigl(\gamma_2 + \frac{\alpha}{\varepsilon}\Bigr) - \cos(\varphi)   \geq \left( \frac{1}{2} - (2-\alpha)(1-\mu)\right)^{1/2}.$$
This inequality implies there cannot be too few points inside the right cone for otherwise $r$ could not attain the size it does. More precisely, this forces
\begin{equation}\label{eq:lowerboundgamma2} \gamma_2 \geq \frac{\cos(\varphi) +  \left( \frac{1}{2} - (2-\alpha)(1-\mu)\right)^{1/2}}{1 + \cos{(\varphi)}} - \frac{\alpha}{\varepsilon}.
\end{equation}

\textbf{Step 4: Using the Hessian.} So far, we have obtained fairly precise information about the number of good points and their approximate location. We have not yet made strong use of the fact that the configuration is a local maximum (we did use it implicitly
when appealing to results of Ling-Xu-Bandeira). In this section, we will explicitly use the fact that the Hessian has to be negative-semidefinite in a maximum to derive a criterion that has to be satisfied for all local maxima. We will then contrast this criterion with the precise structure we have derived to conclude that for some parameters this condition is violated thus excluding them. This will prove the result.
We conclude the argument by using that the configuration of points considered in Step 3 is a 
local maximum: this means that the Hessian is definite which in our
case implies that for any $(w_1, \dots, w_n) \in \mathbb{R}^n$, we
have
\begin{equation} \label{positivity}
 \sum_{i,j}{ a_{ij} \cos{(\theta_i - \theta_j)} (w_i - w_j)^2} \geq 0.
 \end{equation}
 We are free to choose the constants $w_i$ and will choose them in a way that makes the quadratic form as small as possible: this means we want to pick $w_i, w_j$ to be quite different when $i,j$ are on different sides of the cone.
We pick these numbers as follows: for a constant $v \in \mathbb{R}$ to be determined
$$ w = \begin{cases} 1 &\qquad \mbox{for the ones on the left cone}\\
v  &\qquad \mbox{for the outliers}\\
-1 &\qquad \mbox{for the ones on the right cone}.
\end{cases}$$
We will bound the expression from above using the information we have about $\gamma_1$ and $\gamma_2$:  we know that $\gamma_1 + \gamma_2$ is not too small, that $\gamma_1 \leq \gamma_2$ and we have a lower bound on $\gamma_2$.
The quadratic form, which we know to be positive, is bounded from above
\begin{align*}
 \frac{1}{2} \sum_{i,j}{ a_{ij} \cos{(\theta_i - \theta_j)} (w_i - w_j)^2} &=  
  \sum_{i~\tiny \mbox{left}}~  \sum_{j~ \tiny \mbox{right}}{ a_{ij} \cos{(\theta_i - \theta_j)} (w_i - w_j)^2}\\
  &+  \sum_{i~\tiny \mbox{outlier}} ~ \sum_{j~\tiny \mbox{right/left cone}}{ a_{ij} \cos{(\theta_i - \theta_j)} (w_i - w_j)^2}
    \end{align*}
    Now, since $|\sin(\varphi)|^2 \leq 1/2$ by assumption, we have
    that the cosine is negative for any pair of points where one is
    contained in the left cone and one is contained in the right
    cone. Indeed, we see that for any pair $i,j$ of good points on different
    sides of the cone, we have
    $$ \cos{(\theta_i - \theta_j)} \leq \cos{(\pi - (2\pi))},$$
    where $\phi$ is the angle introduced above.
We further bound the quadratic form from above by assuming
    that the number of connections running across good points on different
    sides of the cone is minimized which leads to an upper bound since
    we know that each such connection contributed a nonpositive number.
    To simplify this step in the argument, we will make one more assumption:
  $$ \gamma_2 \geq 1 - \mu.$$ 
Using this assumption, it becomes possible to determine the minimal configuration: each good point on the left-hand side
is connected to all points except the $\gamma_2$ points in the right cone as much as possible. Of course, once $\gamma_2 > 1-\mu$,
each good point on the left has to connect to at least $(\gamma_2 - (1-\mu))n$ good points on the right side of the cone since a point can only be not connected to at most $(1-\mu)n$ other points. This results in
  \begin{align} \sum_{i~\tiny \mbox{left}}  ~~\sum_{j~\tiny \mbox{right}}{ a_{ij} \cos{(\theta_i - \theta_j)} (w_i - w_j)^2} & \leq \gamma_1 (\gamma_2 - (1-\mu)) \cos{(\pi - 2\varphi)} (1 -(-1))^2 n^2 \notag \\
                                                                                                                             & = 4 \gamma_1 (\gamma_2 - (1-\mu)) \cos{(\pi - 2\varphi)} n^2. \label{eq:sumleftright}
  \end{align}
  As for outliers, we have no control on where they are. Let $i$ be an
  outlier and consider the quantity that we need to bound
  $$ S =    \sum_{j~\tiny \mbox{left/right cone}}{ a_{ij} \cos{(\theta_i - \theta_j)} (w_i - w_j)^2}.$$
We can increase this contribution by assuming that the outlier is somewhere in $\varphi \leq \theta_i \leq \pi - \varphi$ and that
all the good points in the left and right cone are located at angles $\varphi$ and $\pi - \varphi$: this means that the good points in the cone are as close to each other as they are allowed to be from the cone condition which simplifies getting large interactions with both of them from the monotonicity of the cosine in that range. If $i$ is an outlier located at angle $\theta$, this leads to the upper bound
$$ S \leq \gamma_1 n (1-v)^2 \cos{(\theta - (\pi - \varphi))} + \gamma_2n (1+v)^2 \cos{(\theta- \varphi)}.$$
We bound expressions of this type via
\begin{align*}
  A \cos{(\theta - (\pi - \phi))} + B \cos{(\theta- \varphi)} &=  -A \cos{(\theta + \varphi)} + B \cos{(\theta - \varphi)} \\
  &= \mbox{Re} \Bigl(-A e^{i (\theta + \varphi)} + B e^{i(\theta - \varphi)}\Bigr)\\
  &= \mbox{Re}~ e^{i \theta} \left( -A  e^{i   \varphi} + B  e^{-i  \varphi} \right)\\
  &\leq  \left|  -A  e^{i   \varphi} + B  e^{-i  \varphi} \right|. \\
  &=   \left|  -A  e^{ 2i   \varphi} + B \right|.
  \end{align*}  
We note that, by assumption, $|\sin{\varphi}|^2 < 1/2$ and thus, since $A,B$ are positive reals,
$$   \left|  -A  e^{ 2i   \varphi} + B \right| \leq \sqrt{A^2+B^2}.$$
This allows us to bound
$$ S \leq n\sqrt{ \gamma_1^2 (1-v)^4 + \gamma_2^2 (1+v)^4}.$$
We now choose the value of $v$ that minimizes this expression. This value is
$$ v = \frac{\gamma_1 - \gamma_2}{\gamma_1 + \gamma_2}$$
and we obtain
$$ S \leq n \frac{4 \gamma_1 \gamma_2 (\gamma_1^2 + \gamma_2^2)^{1/2}}{(\gamma_1 + \gamma_2)^2} = n \frac{4\gamma_1 \gamma_2}{\gamma_1 + \gamma_2} \frac{\sqrt{\gamma_1^2 + \gamma_2^2}}{\gamma_1 + \gamma_2} \leq n \frac{4\gamma_1 \gamma_2}{\gamma_1 + \gamma_2}. $$
Altogether, summing over all the outliers  shows 
$$  \sum_{i~\tiny \mbox{outlier}}  \sum_{j~\tiny \mbox{left/right}}{ a_{ij} \cos{(\theta_i - \theta_j)} (w_i - w_j)^2} \leq \frac{\alpha}{\varepsilon} 
   \frac{4\gamma_1 \gamma_2}{\gamma_1 + \gamma_2}   n^2.$$
  We note that  \eqref{eq:gamma12} implies an upper bound on the number of outliers and
  $$ \frac{1}{\varepsilon} \frac{1}{\gamma_1 + \gamma_2} \leq \frac{1}{\varepsilon}\frac{1}{1 - \frac{\alpha}{\varepsilon}} = \frac{1}{\varepsilon - \alpha}.$$
Combining this with \eqref{eq:sumleftright}, we reach a contradiction to \eqref{positivity} if
 \begin{equation} \label{condition}
    \gamma_1(\gamma_2 - (1-\mu)) \cos{(\pi - 2\varphi)}  + \frac{\gamma_1 \gamma_2  \alpha}{\varepsilon- \alpha}  < 0.
    \end{equation}
 Indeed, if this inequality is satisfied, then the quadratic form corresponding to the Hessian is not definite implying that the configuration we are in does not correspond to a local maximum. We will now analyze the condition.\\
 
 \textbf{Step 5: Analyzing the Condition.} We will now try to understand under what conditions \eqref{condition} holds.
 It clearly requires $\gamma_1 > 0$. We first show that if $\gamma_1 = 0$, then the only stable configuration that can arise
 is actually the one where all points are in the same spot. Then we deal with the more elaborate case that arises when $\gamma_1 > 0$. \\
 
 \textit{Part 1: $\gamma_1 > 0$.} We will show that if $\gamma_1 = 0$, then $\|r\|$ has to be quite big and this will allow us to immediately deduce the desired statement via the Ling-Xu-Bandeira framework. We now discuss this in greater detail.
 If $\gamma_1 = 0$, then there are at least $(1- \frac{\alpha}{\varepsilon})n$ points in the right cone since all the good points are in one of the two cones, none of them in the left cone and there are at least that many good points in total.
 Since the opening angle is
 less than $45^{\circ}$ (recall that this is an assumption and will be the case for all the parameter we will consider below), we know that the $x-$coordinate of $e^{i \theta_j}$ for each good point is at least $1/\sqrt{2}$.  The $x-$coordinate of an outlier is, trivially, at most $-1$. Therefore
 $$ \|r\| = \left|  \sum_{j=1}^{n}{ e^{i \theta_j}} \right| \geq \mbox{Re}   \sum_{j=1}^{n}{ e^{i \theta_j}}  \geq  \frac{1}{\sqrt{2}} \left( 1 - \frac{\alpha}{\varepsilon}\right)n - \frac{\alpha}{\varepsilon}n.$$
 Recalling \eqref{bound}, we have, for any $i$, that
$$  \left| \sin{(\theta_i - \theta_r)} \right| \leq \frac{(1-\mu)n}{\|r\|} \leq \sqrt{2} \frac{1-\mu}{1 - (1+\sqrt{2})\frac{\alpha}{\varepsilon}}.$$
Recalling our regime of interest, $0 \leq \alpha \leq 0.0537$ and $\mu \geq 0.78$ as well as our parameter selection $\varepsilon =0.5$,
we see that
$$  \left| \sin{(\theta_i - \theta_r)} \right| \leq \frac{1}{\sqrt{2}} \frac{1-\mu}{1 - 2\frac{\alpha}{\varepsilon}} \leq 0.45 < \frac{1}{\sqrt{2}}.$$
This case thus reduces to the Proposition of Ling-Xu-Bandeira discussed above and we see that the only possible case is that  all the points are in the same spot.\\
 
\textit{Part 2: Using $\gamma_1 > 0$.} We can thus assume that $\gamma_1 > 0$ in which case the configuration is clearly not the one where all the points are in one spot. We obtain a contradiction to \eqref{positivity} if
 $$   (\gamma_2 - (1-\mu)) \cos{(\pi - 2\varphi)}  + \frac{ \gamma_2  \alpha}{\varepsilon- \alpha}  < 0.$$
Note that $\cos{(\pi - 2\varphi)} = - \cos(2 \varphi)  =  -1 + 2 \sin(\varphi)^2$.
%  we can thus rewrite the condition in terms of  $\sin(\varphi)^2$. 
By assumption, we know that $\sin(\varphi)^2 < 1/2$ and thus the
cosine contribution in the above inequality is negative.  We reach a
contradiction if
  $$ \gamma_2 \left(  \cos{(\pi - 2\varphi)}  + \frac{\alpha}{\varepsilon - \alpha} \right) < (1-\mu)  \cos{(\pi - 2\varphi)}.$$
  We would like to extract a bound for $\gamma_2$ from this but we cannot simply divide by a term without knowing its
  sign which we now determine. 
  The right-hand side is negative; this means that in order to reach a
  contradiction, we certainly would need to require the quantity in
  the parentheses to be negative, i.e.
\begin{equation} \label{condition2}
   \cos{(\pi - 2\varphi)}  + \frac{\alpha}{\varepsilon - \alpha}  < 0.
   \end{equation}
Once this is the case, we can divide by that term and deduce that we reach a global contradiction if
\begin{equation} \label{condition3}
 \gamma_2 \geq \frac{ (1-\mu)  \cos{(\pi - 2\varphi)} }{  \cos{(\pi - 2\varphi)}  + \frac{\alpha}{\varepsilon - \alpha} }.
 \end{equation} 
  Notice that this is condition implies the weaker condition
  $\gamma_2 \geq (1 - \mu)$ that we assumed above, and thus the latter will be
  dropped. Summarizing, we have derived a lower bound on $\gamma_2$ that leads to a global contradiction. At this point, we recall that we did already derive a lower bound on $\gamma_2$ in
  \eqref{eq:lowerboundgamma2} stating that
\begin{align*}
\gamma_2 & \geq \frac{\cos(\varphi) +  \left( \frac{1}{2} - (2-\alpha)(1-\mu)\right)^{1/2}}{1 + \cos{(\varphi)}} - \frac{\alpha}{\varepsilon} \\
& = 1 - \frac{\alpha}{\varepsilon} -  \frac{1 - \left( \frac{1}{2} - (2-\alpha)(1-\mu)\right)^{1/2}}{1 + \cos{(\varphi)}}.
\end{align*}
It is not clear that if our already established lower bound \eqref{eq:lowerboundgamma2} is larger than the lower bound leading to a contradiction, i.e. if
$$  1 - \frac{\alpha}{\varepsilon} -  \frac{1 - \left( \frac{1}{2} - (2-\alpha)(1-\mu)\right)^{1/2}}{1 + \cos{(\varphi)}}  
 \geq \frac{ (1-\mu)  \cos{(\pi - 2\varphi)} }{  \cos{(\pi - 2\varphi)}  + \frac{\alpha}{\varepsilon - \alpha} },
$$
then we have obtained a contradiction. Any collection of points with these parameters must necessarily give rise to a Hessian with a negative definite eigenvalue and thus cannot be a local maximum.\\
 
 \textbf{Summary.}
 In order to obtain a contradiction, it suffices to find, for each $\alpha \leq 0.0537$ two variables $\varepsilon$ and $\delta$
 $$ \alpha < \varepsilon < \delta < 1$$
 such that, abbreviating once again,
 $$ s_{\delta} =  \frac{1}{\sqrt{2}} \sqrt{ \sqrt{9 - 4 \delta} - 3 + 2\delta},$$
 the following properties hold
  \begin{enumerate}
 \item we have
 $$\frac{ \left( s_{\delta} + \frac{\varepsilon}{\delta} (1-s_{\delta})\right)^2 (1-\mu)^2}{  \left( \frac{1}{2} - (2-\alpha)(1-\mu)\right)} < \frac{1}{2}.$$
 This says that the parameters define a critical angle $\varphi = \varphi(\alpha, \mu, \varepsilon, \delta)$ corresponding to an angle less than $45^{\circ}$ which is required for our argument.

 \item we require that this angle  satisfies
   $$   \cos{(\pi - 2\varphi)}  + \frac{\alpha}{\varepsilon - \alpha}  < 0$$
which was required when dividing by it and flipping the sign in Step 5, Part 2.
 \item we also require the angle $\varphi$ satisfies
    $$ 1 - \frac{\alpha}{\varepsilon} - \frac{1 - \left( \frac{1}{2} - (2-\alpha)(1-\mu)\right)^{1/2}}{1 + \cos{(\varphi)}} > \frac{ (1-\mu)  \cos{(\pi - 2\varphi)} }{  \cos{(\pi - 2\varphi)}  + \frac{\alpha}{\varepsilon - \alpha} }.$$
    This shows that the lower bound \eqref{eq:lowerboundgamma2} we derived for $\gamma_2$ is big enough to imply a contradiction, we refer to Step 5, Part 2.\\
 \end{enumerate}
 
 \textbf{Conclusion.} We distinguish the cases $\alpha \leq 0.0537$ and $\alpha > 0.0537$. If $\alpha > 0.0537$, then we immediately obtain a contradiction if
$$ \mu \geq 0.788897.$$
This follows from rerunning the Ling-Xu-Bandeira argument, described in \S 2.1. verbatim and using the definition of $\alpha$, and the value of $\mu$ comes from \eqref{length} and \eqref{bound}.
Let us now assume that $\alpha \leq 0.0537$.
We set
$$ \varepsilon = 0.5 \qquad \mbox{and} \qquad \delta = 0.88.$$
An easy (mathematica) check shows that we obtain a contradiction for the entire range  $0 \leq \alpha \leq 0.0537$ and $0.788897 \leq \mu \leq 0.794$, where $0.794$ is the bound proved in \cite{ling}. More precisely, the inequalities are true with room to spare: we have, over this entire parameter range 
 $$\frac{ \left( s_{\delta} + \frac{\varepsilon}{\delta} (1-s_{\delta})\right)^2 (1-\mu)^2}{  \left( \frac{1}{2} - (2-\alpha)(1-\mu)\right)}  < 0.46 < \frac{1}{2}$$
    $$   \cos{(\pi - 2\varphi)}  + \frac{\alpha}{\varepsilon - \alpha}  < -0.05 < 0$$
and
    $$ 1 - \frac{\alpha}{\varepsilon} - \frac{1 - \left( \frac{1}{2} - (2-\alpha)(1-\mu)\right)^{1/2}}{1 + \cos{(\varphi)}} > \frac{ (1-\mu)  \cos{(\pi - 2\varphi)} }{  \cos{(\pi - 2\varphi)}  + \frac{\alpha}{\varepsilon - \alpha} } + 0.004 .$$
    It is the last expression which is almost satisfied and does not allow an extension to smaller parameter ranges (the critical values occur when $\alpha \sim 0.0537$ and $\mu \sim 0.7889$), the inequality is satisfied with a much bigger gap away from these parameters.
   \end{proof}

\textbf{Acknowledgment.} We are grateful to the anonymous referees for the substantial and detailed reports resulting in a greatly improved manuscript.

\end{document}